\documentclass[11pt]{article}
\usepackage{amssymb}
\usepackage{amsmath}
\begin{document}
\newtheorem{defn0}{Definition}[section]
\newtheorem{prop0}[defn0]{Proposition}
\newtheorem{thm0}[defn0]{Theorem}
\newtheorem{lemma0}[defn0]{Lemma}
\newtheorem{coro0}[defn0]{Corollary}
\newtheorem{exa}[defn0]{Example}
\newtheorem{exe}[defn0]{Exercise}
\newtheorem{rem0}[defn0]{Remark}
\numberwithin{equation}{section}
\def\rig#1{\smash{ \mathop{\longrightarrow}\limits^{#1}}}
\def\swar#1{\swarrow\rlap{$\vcenter{\hbox{$\scriptstyle#1$}}$}}
\def\lswar#1{\swarrow\llap{$\vcenter{\hbox{$\scriptstyle#1$}}$}}
\def\sear#1{\searrow\rlap{$\vcenter{\hbox{$\scriptstyle#1$}}$}}
\def\lsear#1{\searrow\llap{$\vcenter{\hbox{$\scriptstyle#1$}}$}}
\def\near#1{\nearrow\rlap{$\vcenter{\hbox{$\scriptstyle#1$}}$}}
\def\dow#1{\Big\downarrow\rlap{$\vcenter{\hbox{$\scriptstyle#1$}}$}}
\def\ldow#1{\Big\downarrow\llap{$\vcenter{\hbox{$\scriptstyle#1$}}$}}
\def\up#1{\Big\uparrow\rlap{$\vcenter{\hbox{$\scriptstyle#1$}}$}}
\def\lef#1{\smash{ \mathop{\longleftarrow}\limits^{#1}}}
\newcommand{\defref}[1]{Def.~\ref{#1}}
\newcommand{\rank}{\operatorname{rank}}
\newcommand{\propref}[1]{Prop.~\ref{#1}}
\newcommand{\thmref}[1]{Thm.~\ref{#1}}
\newcommand{\lemref}[1]{Lemma~\ref{#1}}
\newcommand{\corref}[1]{Cor.~\ref{#1}}
\newcommand{\exref}[1]{Example~\ref{#1}}
\newcommand{\secref}[1]{Section~\ref{#1}}
\newcommand{\qedd}{\hfill\framebox[2mm]{\ }\medskip}
\def\P#1{{\bf P}^#1}
\def\I{{\cal I}}
\def\O{{\cal O}}
\def\C{{\bf C}}
\def\R{{\bf R}}
\def\Z{{\bf Z}}
\def\proof{{\it Proof:\  \ }}
\def\M{{\cal M}}
\def\F{{\cal F}}
\def\Q{{\bf Q}}
\def\N{{\bf N}}
\date{\ }
\author{ Giorgio Ottaviani}
\title{Symplectic bundles on the plane, secant varieties and L\"uroth quartics revisited}
\maketitle
\vspace*{-0.55in}
\begin{abstract}

\noindent 
Let $X={\bf P}^2\times{\bf P}^{n-1}$ embedded with $\O(1,2)$. We prove that its $(n+1)$-secant variety
$\sigma_{n+1}(X)$ is a hypersurface, while it is expected that it fills the ambient space.
The equation of $\sigma_{n+1}(X)$ is the symmetric analog of the Strassen equation.
When $n=4$ the determinantal map takes $\sigma_5(X)$ to the hypersurface of L\"uroth quartics,
which is the image of the Barth map studied by LePotier and Tikhomirov. 
This hint allows to obtain some results on the jumping lines and the Brill-Noether loci
of symplectic bundles on ${\bf P}^2$
by using the  higher secant varieties of $X$.

\vskip 3pt

\end{abstract}

\section{Introduction}

A L\"uroth quartic is a quartic plane curve which has a inscribed complete pentagon.
By a naive dimensional count (see the Rem. \ref{lurcontracted}), one expects that the general quartic curve has such a inscribed pentagon.
The classical theorem of L\"uroth, published in 1869, says that 
a L\"uroth quartic has not only one, but infinitely many inscribed pentagons.
In equivalent way, L\"uroth quartics form a hypersurface in the space of quartics, so destroying
the naive numerical expectation.
One century later, Barth showed \cite{Bar} the remarkable result that the curve of jumping lines of a rank 2 stable bundle
with $(c_1,c_2)=(0,4)$ on ${\bf P}^2$ is a L\"uroth quartic.   L\"uroth's theorem became again popular,
and Barth gave in \cite{Bar} a new proof of it by using vector bundles.
The equation of the L\"uroth hypersurface is a interesting $SL(3)$-invariant,
and its degree, which is $54$, 
 turned out to be the first nontrivial Donaldson invariant of the plane.

Denote by $M(2,n)$ the moduli space of stable $2$-bundles on $\P 2$ with  $(c_1,c_2)=(0,n)$.
The morphism which associates to any $E\in M(2,n)$ its curve of jumping lines is called the Barth map.

Le Potier and Tikhomirov finally showed \cite{PT} that, conversely, a general L\"uroth quartic is the curve of jumping lines
of a unique bundle $E\in M(2,4)$.  
 They generalized this fact by showing that the Barth map from $M(2,n)$ to the variety of plane curves of degree $n$
 is generically injective.  
\medskip

The second subject of this paper is the higher secant varieties.
For a projective variety $W\subset\P m$ of dimension $n$, the $k$-secant variety $\sigma_k(W)$
is the closure of the locus spanned by $k$ independent points in $W$. 
The expected dimension of $\sigma_k(W)$ is $\min\{m,kn+(k-1)\}$. When $\dim \sigma_k(W)<\min\{m,kn+(k-1)\}$
we say that $W$ is $k$-defective. 
 
A theorem due to Strassen \cite{S} says that $Y=\P 2\times\P {{n-1}}\times\P {{n-1}}$ with its Segre embedding
is $k$-defective for $n$ odd, $n\ge 3$, $k=\frac{3n-1}{2}$ . Indeed, despite of the numerical naive expectation
that $\sigma_k(X)$ fills the ambient space, Strassen proved that  $\sigma_k(X)$ is a hypersurface
and gave an equation for it, see the Rem. \ref{strasseq}.  This equation is called the Strassen equation and it
was later generalized in the papers \cite{LM} and
\cite{LW},
where it was put in the setting of invariant theory.
In \cite{CGG1} and in \cite{AOP} it was studied the defectivity of Segre varieties, and 
in \cite{AOP} it was proposed a conjecture
about which Segre varieties are indeed defective. The equations of
higher secants of Segre varieties are studied, among others, in \cite{LM} and \cite{LW}. The case 
$\P 2\times\P {{n-1}}\times\P {{n-1}}$
 for $n$ odd is probably the most interesting known class of defective Segre varieties,
and its description escaped the geometric techniques working in the other cases.

\medskip

The starting point of this paper was the observation that Strassen equations
are (almost) identical to the equations called $\alpha 3$ by Barth in \cite{Bar}. 
This is a interesting link between secant varieties and vector bundles,
and we tried to explore further this connection. 

A linear algebra approach to $Y=\P 2\times\P {{n-1}}\times\P {{n-1}}$ in the spirit of \cite{Bar}
gives quickly a proof that Strassen equation define the hypersurface $\sigma_k(Y)$
for $n$ odd, $n\ge 3$, $k=\frac{3n-1}{2}$. This is explained in Section 3 (see \thmref{strassen}).
The main argument is the following. Let $\P 2\times\P {{n-1}}\times\P {{n-1}}=\P{{}}(U)\times \P{{}} (V)\times \P{{}} (W)$.
For any $\phi\in U\otimes V\otimes W$, we consider the contraction operator

$$A_{\phi}\colon U\otimes V^{\vee}\rig{}\wedge^2 U\otimes W$$

If  $\phi\in Y$ then $rk A_{\phi}=2$, it follows that if  $\phi\in\sigma_k(Y)$ then $rk A_{\phi}\le 2k$.
So, for $k=\frac{3n-1}{2}$, the tensors $\phi\in\sigma_k(Y)$ have a degenerate $A_{\phi}$
and cannot be general. The Strassen equation is simply
$$\det A_{\phi}=0$$

A natural generalization to the symmetric case is given in Section 4.

The symmetric case corresponds to the Segre-Veronese variety $X=\P 2\times\P {{n-1}}$ embedded by
$\O (1,2)$. We prove that it is $k$-defective for $n\ge 4$,
$n$ even, $k=\frac{3n}{2}-1$,
see \thmref{simstrassen}.

Note that the defective cases appear with $n$ odd in the Segre case and with $n$ even in the 
Segre-Veronese case, adding a touch of intrigue.

Our next result is a proof of L\"uroth's theorem (\thmref{n4})
by using the defectivity of the above Segre-Veronese $X$. In other words, the failure of the numerical expectation
in L\"uroth theorem and in Strassen theorem are two faces of the same phenomenon.

\medskip

The vector bundles $E\in M(2,n)$ were studied in \cite{Bar} as
cohomology bundles of Barth monads
$$I^{\vee}\otimes \O\rig{}V^{\vee}\otimes\Omega^1(2)\rig{f}V\otimes \O(1) $$

where $V$ is a vector space of dimension $n$   and $I$ is a vector space of dimension
$n-2$. 
The L\"uroth quartics are defined by the symmetric determinantal morphism $\Delta$
of $f\in \P{{}}(U\otimes S^2V)$ when $n=4$ (see \cite{D} and \cite{DK}). 
A easy but crucial remark is that the linear map $H^0(f)$ is the symmetric analog of $A_f$,
and it is denoted by $S_f$. This allows to deepen the link between vector bundles
on $\P 2$ and secant varieties to $X$.

A plane curve of degree $k$ which is circumscribed to a complete $(k+1)$-gon is called a Darboux curve. 
It is easy to show, by using the theory of secants, that there are only finitely many 
inscribed complete $(k+1)$-gon in a Darboux curve.
We point out how a question
of Ellingsrud, Le Potier, Stromme and Tikhomirov about the uniqueness of these inscribed $(n+1)$-gons
 is related to the non weak $(n+1)$-defectivity of the Segre-Veronese variety $\P 2\times\P {{n-1}}$
embedded with $\O (1,2)$, (see \thmref{reform}) by a result of Chiantini and Ciliberto, \cite{CC2} .

The Darboux curves correspond, through the Barth map, to $E\in M(2,n)$ such that $E(1)$ has at least one section.
We will show (\thmref{hulsbergen}) that they are defined through the determinantal morphism $\Delta$ from
$\sigma_{n+1}(X)$, so generalizing the L\"uroth case.

The natural objects that come from this description are the symplectic
bundles on $\P 2$. In the Section $7$ we construct the moduli space
$M_{sp}(r,n)$ of symplectic bundles on $\P 2$ of rank $r$ with $c_2=n$ and prove that it is irreducible.
This result follows the lines of \cite{Hu} and it is certainly well known to experts, but we do not know a reference for it.

In the Section $8$ we study the Barth map in the higher rank case,  and we obtain some numerical results on 
it and on Brill-Noether loci in  $M_{sp}(r,n)$. When  $r=2$ we get
an alternative approach to the results of \cite{Bar}.
We describe the curves of jumping lines of a general $E\in M_{sp}(n-2,n)$, see the \thmref{jumpn-2}.

This note has the purpose to serve as a introduction to the subject. So
we tried to give elementary proofs of several of the results that we use, even in the case where they follow from more
 advanced results available in the literature.   An exception to this phylosophy is the Beilinson Theorem,
which we believe is a basic tool, so that it is useful to practice it from the beginning.

I thank Maria Virginia Catalisano for helpful discussions concerning the defectivity of Segre-Veronese varieties and
Pietro Pirola who provided the examples in \propref{pirola}.

\section{Notations and generalities on higher secant varieties}

Let $V$ be a complex vector space.
We denote by $\P{{}}(V)$ the projective space of lines in $V$,
so that $H^0(\P{{}}(V), \O(1))=V^{\vee}$.

For every $a\in\R$, denote by $\lceil a\rceil$ the smallest integer greater or equal than $a$.

We recall that for a projective variety $W\subset\P m$, the $k$-secant variety $\sigma_k(W)$
is the Zariski closure of the set $\bigcup_{x_1,\ldots ,x_k\in W}<x_1,\ldots ,x_k>$,
where $<x_1,\ldots ,x_k>$ denotes the linear span of $x_1,\ldots ,x_k$.

A space spanned by $x_1,\ldots ,x_k\in W$ is called a $k$-secant space.
The expected dimension of $\sigma_k(W)$ is $\min\{m,kn+(k-1)\}$ and it always holds the inequality
$\dim \sigma_k(W)\le\min\{m,kn+(k-1)\}$. When $\dim \sigma_k(W)<\min\{m,kn+(k-1)\}$
we say that $W$ is $k$-defective. $W$ is called defective if it is $k$-defective for some $k$.
 
We can define the abstract secant variety $\sigma^k(W)$
as the Zariski closure of the incidence variety
$$\sigma^k_0(W)=\{(p,x_1,\ldots ,x_k)\in \P m\times Sym^kW | p\in <x_1,\ldots ,x_k>, \dim<x_1,\ldots ,x_k>=k\}$$

It is easy to check (see \cite{Ru} or \cite{Z}) that $\sigma_k(W)$ is the image of the projection on the first factor of
$\sigma^k(W)$.
Since $\sigma^k_0(W)$ is fibered over a open subset of $Sym^kW$ with fibers isomorphic to $\P{{k-1}}$, it follows that $\sigma^k(W)$
is irreducible of dimension $kn+(k-1)$.

\begin{prop0}\label{infty}
If $\dim \sigma_k(W)=kn+(k-1)-d$, then the general point in $\sigma_k(W)$ belongs to $\infty^d$
$k$-secant spaces of dimension $k$.
\end{prop0}

\proof By assumption the general fibers of the projection $\sigma^k(W)\rig{} \sigma_k(W)$
have dimension $d$. \qedd

For improvements of the above proposition see \cite{Ru} or \cite{Z}.

\begin{defn0}\label{wdv}\cite{CC1} $W$ is called $k$-weakly defective if the general hyperplane
which is tangent in $k$ points is tangent along a  variety of positive dimension.
\end{defn0}

The terminology is justified by a result due to Terracini,
who proved that $k$-defective varieties are also $k$-weakly defective (see\cite{CC1}).
In \cite{CC1} are provided counterexamples for the converse and there is a classification of
$k$-defective varieties in small dimension.

\begin{defn0}\label{secord}\cite{CC2}
 The $k$-secant order $d_k(W)$ is  the number of irreducible components of the general fiber of
the morphism $\sigma^k(W)\rig{} \sigma_k(W)$.
\end{defn0}

When  $\dim \sigma_k(W)=kn+(k-1)$, then $d_k(W)$ is the degree of the morphism
$\sigma^k(W)\rig{} \sigma_k(W)$.

$d_k$ measures how many $k$-secant spaces pass through a general point in $\sigma_k(W)$.

The secant order $d_k$ is called $\mu_{k-1}$ in \cite{CC2}, due on a different terminology
about higher secant varieties. We choosed to change the
letter in order to avoid any confusion.

\begin{thm0}[Chiantini-Ciliberto, \cite{CC2}, Corollary 2.7]
\label{maincc}Let $W\subset\P m$ be a variety such that $\dim \sigma_k(W)=kn+(k-1)<m$.
Then $d_k(W)=1$ unless it is $k$-weakly defective.
\end{thm0}

\section{Strassen equations and the tensors $3\times n\times n$}
\label{3nn}

The aim of this section is to give a simpler proof of the defectivity of $\P 2\times\P {{n-1}}\times\P {{n-1}}$
\cite{S}.

Let $U$,$V$, $V'$ be vector spaces of dimension respectively $3$, $n$, $n$. The variety of the decomposable tensors in the projective space
$\P{{}} (U\otimes V\otimes V')$ is the Segre variety
$Y=\P{{}}(U)\times \P{{}}(V)\times \P{{}}(V')$.

The expected minimal $p$ such that $\sigma_p(Y)=\P{m}$
is
$$p=\lceil\frac{3n^2}{2n+1}\rceil$$

It is elementary to check that

\begin{equation}\label{maxp0}
2p=\left\{\begin{array}{ll}3n-1&\hbox{\ if\ }n\hbox{\ is odd}\\
3n&\hbox{\ if\ }n\hbox{\ is even\ }, \\
\end{array}\right.
\end{equation}

For any $\phi\in U\otimes V\otimes V'$, we consider the contraction operator

$$A_{\phi}\colon U\otimes V^{\vee}\rig{}\wedge^2 U\otimes V'$$

If $P$, $Q$, $R$ are the three $n\times n$ slices of $\phi$,
the matrix of $A_{\phi}$ in the obvious coordinate system is

$$\left[\begin{array}{ccc}
0&P&Q\\
-P&0&R\\
-Q&-R&0\\
\end{array}\right]$$

\begin{lemma0}\label{eqz} Let $Q$ be invertible. Then

$$\left[\begin{array}{ccc}
I&0&0\\
0&I&-PQ^{-1}\\
0&0&I\\
\end{array}\right]
\left[\begin{array}{ccc}
0&P&Q\\
-P&0&R\\
-Q&-R&0\\
\end{array}\right]
\left[\begin{array}{ccc}
I&0&0\\
0&I&0\\
0&-Q^{-1}P&I\\
\end{array}\right]=
\left[\begin{array}{ccc}
0&0&Q\\
0&Z&R\\
-Q&-R&0\\
\end{array}\right]
$$

where
$Z=PQ^{-1}R-RQ^{-1}P=Q[Q^{-1}P,Q^{-1}R]$.
\end{lemma0}
\proof Straightforward. \qedd

From the \lemref{eqz} it follows in particular
$$\det\left[\begin{array}{ccc}
0&P&Q\\
-P&0&R\\
-Q&-R&0\\
\end{array}\right]=(\det Q)^2\det\left(PQ^{-1}R-RQ^{-1}P\right)$$

and hence the formula
$$\det\left[\begin{array}{ccc}
0&P&Q\\
-P&0&R\\
-Q&-R&0\\
\end{array}\right](\det Q)^{n-2}=\det\left(P\cdot adj(Q)\cdot R-R\cdot adj(Q)\cdot P\right)$$
which holds for any $P$, $Q$, $R$.

\begin{thm0}
\label{strassen}

(i)  If  $\phi\in\sigma_k(\P 2\times\P {{n-1}}\times\P {{n-1}})$ then $rk A_{\phi}\le 2k$.

(ii) If $\phi$ is generic and $n\ge 3$,  then $rk A_{\phi}={3n}$, hence $\P 2\times\P {{n-1}}\times\P {{n-1}}$ is $k$-defective
when $n$ is odd and $k=\frac{3n-1}{2}$.

\end{thm0}

\proof 

If $\phi$ is decomposable, say $\phi=u_1\otimes v_1\otimes w_1$ then $rk A_{\phi}=2$, indeed $Im\ A_{\phi}$
is $(u_1\wedge U)\otimes w_1$.

It follows that if $\phi\in\sigma_k(Y)$ then $rk A_{\phi}\le 2k$, proving (i).

Let $\lambda_i, \mu_j$ be generic constants. With obvious notations for basis in the vector spaces
$U$, $V$, $V'$, let $\phi=u_1\otimes(\sum \lambda_i v_i\otimes w_i)+u_2\otimes(\sum v_i\otimes w_i)
+u_3\otimes(\sum \mu_i v_i\otimes w_{i+1})$,
where we denote $w_{n+1}=w_1$. Then with the matrix notations of \lemref{eqz} we have
$P=diag(\lambda_i)$, $Q=Id$ and
$$R=\left[\begin{array}{cccc}
&\mu_1\\
&&\ddots\\
&&&\mu_{n-1}\\
\mu_n\\
\end{array}\right]$$
Hence 
$$Z=[P,R]=
\left[\begin{array}{cccc}
&\mu_1(\lambda_1-\lambda_2)\\
&&\ddots\\
&&&\mu_{n-1}(\lambda_{n-1}-\lambda_n)\\
\mu_n(\lambda_n-\lambda_1)\\
\end{array}\right]$$
It follows that $Z$ is invertible and moreover by \lemref{eqz} $rk A_{\phi}={3n}$, proving (ii)
The $k$-defectivity for $n$ odd and $k=\frac{3n-1}{2}$ follows by (\ref{maxp0}).

\qedd

\begin{rem0}\label{strasseq} $\sigma_k(Y)$ is the hypersurface with equation
$\det A_{\phi}=0$.
This is proved in \cite{S} lemma 4.4 by computing the tangent spaces at   $\frac{3n-2}{2}$
suitably chosen points in $Y$. 
By \lemref{eqz} this equation is equivalent to the Strassen one which has the nice commuting shape 
$$\det\left(P\cdot adj(Q)\cdot R-R\cdot adj(Q)\cdot P\right)=0$$
which has to be divided by $(\det Q)^{n-2}$ in order to get a
$SL(U)\times SL(V)\times SL(V')$-invariant form with homogeneous weights.
\end{rem0}
\begin{rem0} Also $\P 2\times\P 3\times\P 3\subset\P {{47}}$ is defective. In \cite{AOP}
it is proved that the dimension of $\sigma_5(\P 2\times\P 3\times\P 3)$ is $43$ instead of $44$
by showing that through five generic points in $\P 2\times\P 3\times\P 3$ there is a rational normal curve
$C_8$. Note that 
the argument in the above proof does not work in this case. 
\end{rem0}
\begin{rem0} Strassen proves in \cite{S} prop. 4.7 that if $n$ is even then
$\sigma_{\frac{3n}{2}}(\P 2\times\P {{n-1}}\times\P {{n-1}})$ is the ambient space by exhibiting
several explicit tensors.
 \end{rem0}
\section{The  symmetric  tensors $3\times n\times n$}
\label{3nns}

The aim of this section is to extend the results of Section \ref{3nn}
to the  symmetric  cases.
The defectivity
of Segre-Veronese varieties has been studied in \cite{CGG2}, \cite{CaCh} and \cite{CarCa}.
Tony Geramita and Enrico Carlini pointed out to me  that the defectivity of $\P 2\times\P {{3}}$
embedded with $\O (1,2)$ was classically known, and it essentially appears in a paper of
Emil  Toeplitz \cite{T}(the father of Otto), 
see \cite{CaCh} theor. 4.3 (and also \cite{CarCa} 6.2), where the reader can find a modern
geometrical proof. 
It turned out that E. Toeplitz  already wrote in 1877 (\cite{T}, pag. 441) the equation
of $\sigma_5(X)$, which is the symmetric analog to the Strassen equation. 
The approach of the previous section allows to generalize this result to $X=\P2\times\P{{n-1}}$.

Let $U$, $V$ be vector spaces of dimension $3$, $n$.

Consider the Segre-Veronese embedding $X=\P{{}}(U)\times\P{{}}(V)\rig{}\P{{}}(U\otimes S^2V)=\P{m}$
with the linear system $\O (1,2)$, so that
$m=3{{n+1}\choose 2}-1$.

The expected minimal $p$ such that $\sigma_p(X)=\P{m}$
is
$$p=\lceil\frac{m+1}{n+2}\rceil$$

It is elementary to check that

\begin{equation}\label{maxp}
2p=\left\{\begin{array}{cl}3n-1&\hbox{\ if\ }n\hbox{\ is odd}\\
3n-2&\hbox{\ if\ }n\hbox{\ is even\ }, n\ge 4\\
6&\hbox{\ if\ }n=2\end{array}\right.
\end{equation}

We can study this embedding in the following way.
Let $\phi\in U\otimes S^2V$. It defines the contraction
$$S_{\phi}\colon U\otimes V^{\vee}\rig{}\wedge^2U\otimes V$$

Note that since $\wedge^2U\simeq U^{\vee}$ for the $SL(U)$-action,
we have
\begin{equation}
\label{sfi}
S_{\phi}^t=-S_{\phi}
\end{equation}

Indeed, if $P$, $Q$, $R$ are the three $n\times n$ symmetric slices of $\phi$,
the matrix of $S_{\phi}$ in the obvious coordinate system is  again

$$\left[\begin{array}{ccc}
0&P&Q\\
-P&0&R\\
-Q&-R&0\\
\end{array}\right]$$

which is now skew-symmetric.
\begin{thm0}
\label{simstrassen}

(i)  If  $\phi\in\sigma_k(X)$ then $rk S_{\phi}\le 2k$.

(ii) If $\phi$ is generic and $n\ge 2$, $n$ is even, then $rk S_{\phi}={3n}$. Hence  we get that $X$ is $k$-defective for $n\ge 4$,
$n$ even, $k=\frac{3n}{2}-1$.

\end{thm0}

\proof 
It is analogous to the proof of \thmref{strassen}.
If $\phi$ is decomposable then $rk S_{\phi}=2$.
The variety of the decomposable tensors in the projective space
$\P{{}} (U\otimes S^2V)$ is the Segre-Veronese variety
$X=\P{{}} (U)\times \P{{}} (V)$.

It follows that if $\phi\in\sigma_k(X)$ then $rk S_{\phi}\le 2k$, proving (i).

Let $n=2h$, $\lambda_i, \mu_j$ be generic constants for $1\le i\le n$, $1\le j\le h$.
 With obvious notations for basis in the vector spaces
$U$, $V$ let $\phi=u_1\otimes(\sum v_i^2)+
u_2\otimes(\sum \lambda_i v_i^2)+u_3\otimes[\sum \mu_i (v_i+v_{i+h})^2]$.
Then $rk S_{\phi}={3n}$, proving (ii) thanks to (\ref{maxp}).
\qedd

{\bf Remark} 
It is surprising that in the \thmref{strassen} and \thmref{simstrassen}
the odd and even cases exchange.
If $n$ is even, $n\ge 4$, the secant variety $\sigma_{\frac{3n-2}{2}}(X)$
is a hypersurface of degree $\frac{3n}{2}$, with equation $Pf \left(S_{\phi}\right)=0$.
This can be proved like in \cite{S} lemma 4.4.

If $n$ is odd then $\sigma_{\frac{3n-1}{2}}(X)$ fills the ambient space, indeed $S_{\phi}$ is always singular.
Following \cite{S} lemma 4.4, the list of tensors such that their tangent spaces span the ambient space is the following.
Let $v_1,\ldots ,v_n$ be a basis of $V$ and  

 $u_i$, $1\le i\le n\qquad$,
$\tilde{u_{\nu}}$, $1\le \nu\le n-2,\quad\nu $ odd
be vectors in $U$ such that their coefficients are algebraically independent over ${\bf Q}$.

Then the list of tensors is

$u_i\otimes v_i^2$, $1\le i\le n$

$\tilde{u_{\nu}}\otimes (v_{\nu}+v_{\nu+1}+v_{\nu+2})^2$, $1\le \nu\le n-2,\quad\nu\ $ odd.

\bigskip

\section{Generalities on plane curves as linear symmetric determinants}
For any $f\in \P{{}}(U\otimes S^2V)$
consider the map ${\tilde f}\colon V^{\vee}\to V\otimes U$
as a map from $V^{\vee}$ to $V$ with coefficients in $U$. It can be represented
as a $d\times d$ matrix with coefficient linear forms on $\P{{}} (U^{\vee})$.

Its determinant $\Delta(f)$  gives a morphism
$$\P{{}}(U\otimes S^2V)\setminus Z(\Delta)\rig{\Delta}\P{{}}(S^nU)$$
where $Z(\Delta)$ is the locus $\{f |\Delta (f)\equiv 0\}$.
The elements of $\P{{}}(S^nU)$ can be regarded as degree $n$ curves in
$\P{{}}(U^{\vee})$.

The morphism $\Delta$ was classically studied as discriminant locus of a linear system of quadrics.
All the coefficients of $\Delta$ are $SL(U)$-invariant. 

Assume now $\dim U=3$, hence $\Delta(f)$, for $f\notin Z(\Delta)$,  is the equation of a plane curve $C$ of degree $n$.
The classical Theorem of Dixon, that we review in this section,
says that the general plane curve $C$ is in the image of $\Delta$, and each symmetric 
determinantal representation of $C$ corresponds to a theta characteristic $\theta$ on $C$ such that
$h^0(\theta)=0$.
In  \cite{Bea1}, prop. 4.2,
the Theorem of Dixon is proved as a consequence of a more general result about ACM bundles.

A general $f\in \P{{}}(U\otimes S^2V)$ gives a subspace $\P{{}}(U^{\vee})\subset\P{{}} (S^2V)$.
The space of quadrics $\P{{}} (S^2V)$ contains the discriminant hypersurface ${\cal D}$ and $C$ is identified with
the intersection
$\P{{}}(U^{\vee})\cap {\cal D}$. Note that $Sing {\cal D}$ consists
of the symmetric matrices of rank $\le n-2$ and it has codimension three. The general immersion
of $\P{{}}(U^{\vee})$ does not meet $Sing {\cal D}$ , hence the curve $C$ is in general smooth
by Bertini Theorem. We denote by $\phi\colon C\to \P{{}} (S^2V)$ the previous immersion,
hence every $x\in C$ defines (up to scalar multiplication) a morphism
$\phi(x)\colon V^{\vee}\to V$, where we used the same letter, with a slight abuse of notation.
This smooth curve of genus $g={{d-1}\choose 2}$ comes equipped with a second immersion in $\P{{}} (V^{\vee})$ that is defined as
$$\begin{array}{cccc}\psi\colon &C&\to &\P{{}} (V^{\vee})\\
&x&\mapsto&\ker \phi(x)\\
\end{array}$$.

We call $\psi^*O(1)=L$. The key result about this second immersion
is that the line bundle $L$ is related to 
 a theta-characteristic of $C$. We follow \cite{D}.
 
 The kernel of $ \phi(x)$ can be computed by taking the adjoint matrix $Ad \phi(x)$.
 Indeed the adjoint of a symmetric matrix of rank $n-1$
 is a symmetric matrix of rank one, hence defining a element in the quadratic Veronese variety.
 This means that the embedding given by $L^2$ is given by the minors of $\phi(x)$ which have degree
 $d-1$. Hence $\deg L^2=d(d-1), \quad L^2=\O_C(d-1)$, and it follows that
 $$\deg L={d\choose 2}$$
 Set $\theta=L(-1)$, we get $$\theta^2=\O_C(d-3)=K_C$$
 that is $\theta$ is a theta-characteristic.

Moreover $L$ is generated by the $d$ sections given by the embedding $\psi$,
that is we have a surjection 
$$\O_{\P{{}}(U^{\vee})}\otimes V^{\vee}\rig{}L\rig{}0$$

We could compute now the kernel of this morphism, but we do not pursue this
because we will follow in a while
another road, by using the Beilinson Theorem. 

Conversely, given a theta-characteristic  $\theta$ on a smooth plane curve $C$ of degree $n$ such that $h^0(\theta)=0$,
we can construct $\theta$ as symmetric linear determinant.

Indeed let $L=\theta(1)$, then we define $V^{\vee}=H^0(L)$ which has dimension $n$ and we get the embedding
$\psi\colon C\to \P{{}} (V^{\vee})$. By composing with the Veronese embedding we get
$C$ in $\P{{}} (S^2V^{\vee})$ with associated line bundle $L^2=K_C(2)=\O(d-1)$.
Since $\O(d-1)$ is the restriction of $\O_{\P{{}}(U)}(d-1)$, it follows that this last embedding is the restriction of a embedding $\P{{}}(U)\subset \P{{}} (S^2V^{\vee})$
given by a linear system of plane curves of degree $d-1$. 
In particular it is defined $f\in \P{{}}(U\otimes S^2V)$. An explicit computation (see \cite{D} pag. 81)
shows that $\Delta(f)$ is the equation of $C$.

The general result is the following
 
\begin{thm0}[Dixon]\label{dixon}
Let $C$ be a smooth plane curve 
defined by a polynomial $F$ of degree $d$ and $\theta$ be a theta-characteristic on $C$ such that $h^0(\theta)=0$.
There is a symmetric map $M$ such that

$$0\rig{}\O(- 2)^d\rig{M}\O(-1)^d\rig{}\theta\rig{}0$$
and $\det~M=F$. Two maps $M$, $M'$ define the same $\theta$ if and only if they lie in the same 
$SL(V)$-orbit.  In particular $\Delta$ is dominant.

Conversely, the cokernel
of a injective symmetric map $\O(- 2)^d\rig{M}\O(-1)^d$
is a theta-characteristic $\theta$ on the curve $C$ defined by $\det M=0$, such that $h^0(\theta)=0$.
\end{thm0}

The quickest proof of the \thmref{dixon} is probably
obtained as an application of the Beilinson Theorem .

We recall the Beilinson theorem in the form given in \cite{AO}.
\begin{thm0}[Beilinson]
\label{beilinson}
Let $F$ be a coherent sheaf on $\P n$ and let $Q$ be the quotient bundle. 
There is a complex $$\ldots\rig{d_{-1}}C^0\rig{d_0}C^{1}\rig{d_1}\ldots$$ on $\P n$ such that

(i) $C^h=\oplus_{j+h=i}\wedge^jQ^{\vee}\otimes H^i(F(-j))$

(ii)
the horizontal maps extracted from $d_i$
$$\wedge^jQ^{\vee}\otimes H^i(F(-j))\rig{}\wedge^{j-1}Q^{\vee}\otimes H^i(F(-j+1))$$
are the natural multiplication maps

(iii) the cohomology is

$$\frac{\textrm{Ker\ }d_h}{\textrm{Im\ }d_{h-1}}=
\left\{\begin{array}{ll}0&\textrm{if\ }h\neq 0 \\
F&\textrm{if\ }h=0\\
\end{array}
\right.
$$
\end{thm0}

{\it Proof of \thmref{dixon}:}
Consider $\theta$ as a coherent sheaf on $\P 2$ extended to zero.

The Beilinson table  for $\theta(1)$

$$\begin{array}{lll}H^2(\theta(-1))&H^2(\theta)&H^2(\theta(1))\\
H^1(\theta(-1))&H^1(\theta)&H^1(\theta(1))\\
H^0(\theta(-1))&H^0(\theta)&H^0(\theta(1))\\
\end{array}$$
is
$$\begin{array}{lll}0&0&0\\
V^{\vee}&0&0\\
0&0&V\\
\end{array}$$
hence we get from \thmref{beilinson} the resolution

$$0\rig{}V^{\vee}\otimes\O(-2)\rig{M}V\otimes\O(-1)\rig{}\theta\rig{}0$$
Applying the functor $\underline{Hom}(-,\O) $ we get
$$ 0\rig{}V^{\vee}\otimes\O(1)\rig{M^t}V\otimes\O(2) \rig{}\underline{Ext}^1(\theta,\O)$$
and by Grothendieck duality ([FGA] pag. 149-08) $$\underline{Ext}^1(\theta,\O)=\theta^{\vee}(3)$$

so that twisting by $\O(-3)$ we get that $M=M^t$, hence the map is symmetric.

A morphism between $\theta$ and $\theta'$ lifts to a morphism between the two resolutions.
The converse follows again by the Grothendieck duality.\qedd

\begin{coro0}\label{dimf} Let $f\in \P{{}}(U\otimes S^2V)$ 
such that $\Delta(f)$ is a smooth plane curve in $\P{{}}(U^{\vee})$. Then the fiber 
of $\Delta$ containing $f$ has dimension
equal to $\dim SL(V)=n^2-1$.
\end{coro0}
\proof By the \thmref{dixon} the fiber is a union of orbits,
hence their dimension is $\le n^2-1$.
Since the map $\Delta$ is dominant every fiber has dimension $\ge n^2-1$. \qedd

{\bf Remark} Moreover {\it every} smooth plane curve is in the image of $\Delta$ (\cite{Bea1}, Remark 4.4)
but we will do not need this fact.

Wall studied in \cite{W} the map $\Delta$ in the setting of invariant theory.
He remarked, as a consequence of Dixon theorem, that the field of invariants
of $SL(V)$ acting on $\P{{}}(U\otimes S^2V)$ is a finite extension
of the field generated by the coefficients of $\Delta$. In other terms, the semistable points
for the action of $SL(V)$ on $\P{{}}(U\otimes S^2V)$
are exactly given by the locus where $\Delta$ is not defined.
Wall also proves that if $\Delta(f)$ is a reduced curve,
then the stabilizer of $f$ is finite.
In particular if $\Delta(f)$ is smooth then the stabilizer of $f$ is finite.
In other terms

\begin{prop0}[Wall]\label{wall} There is a factorization
 through
the GIT quotient
$$\begin{array} {ccc}
\P{{}}(U\otimes S^2V)^{ss}\\
\dow{\pi} &\sear{\Delta}\\
\P{{}}(U\otimes S^2V)//SL(V)&\rig{g}&\P{{}}(S^nU)\\
\end{array}$$

where $g$ is generically finite. 
\end{prop0}

\section{L\"uroth quartics revisited}

\begin{defn0} A complete $n$-gon is the union of $n$ lines
in $\P 2=\P(U^{\vee})$ meeting in ${n\choose 2}$ distinct vertices.
\end{defn0}

\begin{defn0}
A L\"uroth quartic is a smooth quartic which has a inscribed complete pentagon,
that is it contains its ten vertices. More generally
a plane curve of degree $n$ which has a inscribed a complete $(n+1)$-gon is called a Darboux curve.
\end{defn0}

Let us identify, again with a slight abuse of notations,
the Segre Veronese variety $X\simeq \P{{}} (U)\times \P{{}} (V)$ embedded by $\O(1,2)$
with its cone in $U\otimes S^2V$.
As we saw in Section \ref{3nns}, the elements in $U\otimes S^2V$ are stratified by the (border) rank,
denoting by rank one the elements of $X$ 
and by rank $k$ the elements of $\sigma_k(X)\setminus \sigma_{k-1}(X)$.

Note that if $f\in X$ then $\textrm{rank} S_f=2$
so that if $f\in \sigma_k(X)$ then $\textrm{rank} S_f\le 2k$.

So we are exactly in the setting of \thmref{simstrassen}.

\begin{prop0}
\label{nn+1}(i)
If  $f\in\sigma_n(X)\setminus \sigma_{n-1}(X)$  is general then the plane curve
given by $\Delta ( f)=0$  consists of $n$ lines.

(ii) If $f\in\sigma_{n+1}(X)\setminus \sigma_n(X)$ is general then the plane curve
of degree $n$ given by $\Delta ( f)=0$ is a Darboux curve (so for $n=4$ is a L\"uroth quartic).  The divisor given by the 
${{n+1}\choose 2}$ vertices of
the  complete inscribed $(n+1)$-gon has the form $\theta+2H$
where $\theta$ is the theta characteristic defined in \thmref{dixon} and $H$ is the hyperplane divisor.
\end{prop0}
\proof (i) Assume $f=\sum_{i=1}^n u_i\otimes v_i^2 $.
The vectors $u_i\in U$ defines lines $L_i$ in the plane $\P{{}}(U^{\vee})$ by the equation $u_i(-)=0$.
Then the matrix corresponding to $f$ evaluated on the line $L_i$ has rank $\le n-1$ and it is degenerate.
This means that $\cup_i L_i\subseteq \{\Delta(f)=0\}$ and the other inclusion holds by degree reasons.

(ii) Assume $f=\sum_{i=1}^{n+1} u_i\otimes v_i^2 $.
Then 
the matrix corresponding to $f$ evaluated in the point of intersection $L_p\cap L_q=0$
is a symmetric matrix of rank $\le n-1$ (because two summands
vanish), hence it is degenerate and the point belongs to the curve $\{\Delta(f)=0\}$.
Let $D$ be the divisor corresponding to the vertices of the $(n+1)$-gon,
so that $\deg D={{n+1}\choose 2}$. Set $\theta=D-2H$. Since every vertex contains
two edges we get that
$2D=(n+1)H$ and it follows that $2\theta = 2D-4H = (n-3)H=K$, hence  
$\theta$ is a theta-characteristic. 
An explicit computation (see the next example (\ref{smoothn}) shows that $\theta$ is the theta characteristic defined in \thmref{dixon}.
\qedd

It is well known that there exist smooth curves of degree $n$ with a
inscribed complete $(n+1)$-gon. With the determinantal representation their equation can be constructed
from general linear forms $l_1,\ldots ,l_n, l_{n+1}=l$ as follows

\begin{equation}\label{smoothn}
\left|\begin{array}{ccccc}
l_1+l&l&l&\ldots&l\\
l&l_2+l&l&\ldots&l\\
\vdots&\vdots&&&\vdots\\
l&l&\ldots&\ldots&l_n+l\\
\end{array}\right|=\sum_{i=1}^{n+1}l_1l_2\cdots\hat{l_i}\cdots l_{n+1}=0
\end{equation}

The following theorem was proved in \cite[Lemma2]{Bar} in a more general setting.
For the convenience of the reader we repeat here Barth's proof in our case.

\begin{lemma0}\label{quartpent}
Given a general complete $(n+1)$-gon, the space of curves of degree $n$ passing through its ${{n+1}\choose 2}$ vertices
has projective dimension $n$, that is the ${{n+1}\choose 2}$ vertices impose independent conditions to curves of degree $n$.
\end{lemma0}
\proof Let $Z$ be the scheme given by the ${{n+1}\choose 2}$ vertices and let $r$ by a general line.

We have the exact sequence
$$0\rig{}H^0({\cal I}_Z(n-1))\rig{}H^0({\cal I}_Z(n))\rig{}H^0({\cal I}_Z(n)_{|L})$$

The first space is zero dimensional because a curve of degree $n$ through $Z$  vanishes
on all the edges of the $(n+1)$-gon. Since ${\cal I}_Z(n)_{|L}\simeq\O_L(n)$
it  follows $h^0({\cal I}_Z(n))\le n+1$. Moreover $h^0({\cal I}_Z(n))\ge {{n+2}\choose 2}-{{n+1}\choose 2}=n+1$,
as we wanted. \qedd

\begin{thm0} (L\"uroth, 1869) \label{n4} If a plane curve of degree $4$  has a inscribed complete pentagon,
then it has $\infty^1$ inscribed complete pentagons. Equivalently, the (closure of the) locus of
L\"uroth quartics is a hypersurface in $\P {{}}(S^4V)$.
\end{thm0}                                  

\proof By \propref{nn+1} every $f\in\sigma_{5}(X)\setminus \sigma_4(X)$
defines a L\"uroth quartic with equation $\Delta(f)$. By \thmref{simstrassen}
the higher secant variety $\sigma_{5}(X)$ is defective and has dimension one less then expected.
 By \propref{infty} $f$ belongs to $\infty^1$ $5$-secant spaces.
The quartics defined in this way are the image through the determinantal morphism
of the $SL(V)$-invariant variety $\sigma_{5}(X)$, so by \thmref{dixon}
(since there are smooth L\"uroth quartics) they give a irreducible hypersurface in $\P{{}}(S^4U)$. 

This is the core of the argument but it does not yet conclude the proof. In order to prove that the general (and hence any)
L\"uroth quartic has $\infty^1$ inscribed pentagons, we need the following easy dimensional count.

Consider the product $\P{{}}(S^4U)\times\left({\P{{}}}(U)\right)^5$ and the incidence variety $I$
given by the closure of 
$$I^0=\{(f,l_1,\ldots ,l_5)|\ l_i \textrm{\ define a complete pentagon inscribed in the 
smooth\ }f\}$$ By \lemref{quartpent} the general fiber of this projection is a linear
space of projective dimension $4$,
hence $I$ is irreducible and $\dim I=14$.
Clearly the projection of $I$ on $\P{{}}(S^4U)$ gives the closure of the locus of L\"uroth quartics.
Our first argument
 showed that
$I$ contains a dense open subset of pairs given by a L\"uroth quartic and its
inscribed pentagons, hence the projection of $I$ is a hypersurface. This concludes the proof. \qedd

\begin{rem0} The above proof is close to Frahm one (cite{Fr}),
see also the Remark 4.5 in \cite{CaCh}, although in these sources the approach is a bit different.
We hope it is still useful to review this subject in the modern language
and to make direct the link with the higher secant variety of $X$.
 The original proof of L\"uroth was based on the polarity with Clebsch
quartics, and it has been reviewed in \cite{DK}.
\end{rem0}
\begin{rem0} It can happen that $f\notin\sigma_{5}(X)$ but still $\Delta(f)$
is a L\"uroth quartic. Indeed $f\notin\sigma_{5}(X)$ picks only one of the $36$ connected components of a general
fiber of $\Delta$.
\end{rem0}
The hypersurface of L\"uroth quartics was classically studied. F. Morley
(the same of trisector theorem) proved in 1918 (\cite{Mo}) that it has degree $54$,
by the classical Aronhold construction of plane quartics from seven given points.

\begin{thm0}\label{expn+1} Let $X= \P2\times\P{{n-1}}$
embedded with the linear system $\O (1,2)$. For $n\ge 5$ the secant variety
$\sigma_{n+1}(X)$ has the expected dimension $(n+1)^2+n$.
\end{thm0}
\proof For $n=5$ it is enough to pick six generic points , compute their tangent spaces and apply the 
Terracini lemma. For $n\ge 6$ we can use the standard inductive technique, which goes back to Terracini.

Let consider the divisor  $X'= \P2\times\P{{n-2}}\subset X$.
Given $n+1$ $P_i$ points on $X$ let specialize the first $n$ of them on $X'$.
Let $Z=\{P_1,\ldots ,P_{n}\},\qquad Z^2=\{P_1^2,\ldots ,P_{n}^2\}$.

We get the exact sequence

$$0\rig{}I_{Z\cup P_{n+1}^2,X}(1,1)\rig{}I_{Z^2\cup P^2_{n+1},X}(1,2)\rig{}I_{Z^2,X'}(1,2)$$

The first space has always the expected dimension. The last space has the expected dimension
by  the inductive assumption. Hence also the middle space has the expected dimension.\qedd

The following theorem says that the statement of the L\"uroth \thmref{n4} is false
for $n\ge 5$.

\begin{coro0} \label{n5}Let $n\ge 5$. The generic Darboux curve of degree $n$ 
has only finitely many inscribed complete $(n+1)$-gons.
\end{coro0}

\proof Let $M^{sm}_n\subset\P{{}}(S^nU)=M_n$ be the variety of smooth plane curves of degree $n$,
 let $G_n$ be the open (smooth) part of $Sym^{n+1}\P2$
parametrizing complete $(n+1)$-gons.
Let $R_n$ be the closure of the incidence variety in $M_n\times G_n$, namely 
\begin{equation}\label{defw}
R_n=\overline{R_n^0}\qquad
R_n^0=\{(m,g)\in M^{sm}_n\times G_n |\ g\hbox{\ is inscribed to\ } m\}
\end{equation}
Let $p$, $q$ be the restriction to $R_n$ of the two projections respectively on $M_n$ and $G_n$.
By \propref{quartpent} the generic fiber of $q$ is isomorphic to $\P n$, hence $R_n$ is irreducible and $\dim R_n=3n+2$.
From \thmref{expn+1}, the examples (\ref{smoothn}) and the \corref{dimf} the image $p(R_n)$ has dimension $\ge (n+1)^2+n-(n^2-1)=3n+2$.
It follows that $p$ is generically finite, as we wanted.
\qedd

\begin{rem0}\label{lurcontracted} The first part of the above proof still says that $R_4$ is irreducible of dimension $14$.
But in this case the \thmref{expn+1} fails and $R_4$ is contracted by $p$ to the L\"uroth hypersurface.
We remark also that a open subset of $R_n$ comes from determinantal curves in $\sigma_{n+1}(X)$
by the construction in the examples (\ref{smoothn}).
\end{rem0}
It follows that

\begin{prop0}\label{darblocus} Let $n\ge 4$.
$\Delta(\sigma_{n+1}(X))$ is the Darboux locus, that is the closure of the
variety consisting of Darboux plane curves of degree $n$. Its dimension is $13$ for $n=4$ and
$3n+2$ for $n\ge 5$.
\end{prop0}

Denote by $\pi_1\colon X\to\P 2$ the projection on the first factor.
The abstract secant variety $\sigma^{n+1}(X)$ has a dominant morphism
$$\Delta\times \pi_1^{n+1}\colon \sigma^{n+1}(X)\to R_n$$

which identifies with the $SL(V)$-quotient $\sigma^{n+1}(X)//SL(V)$.
To prove this claim, let consider, like in \propref{wall}, the factorization

$$\begin{array} {ccc}
\sigma^{n+1}(X)^{ss}\\
\dow{\pi} &\sear{\Delta\times \pi_1^{n+1}}\\
\sigma^{n+1}(X)//SL(V)&\rig{g}&R_n\\
\end{array}$$

where $g$ is finite. Since the inscribed $(n+1)$-gon determines 
by \propref{nn+1} the theta divisor,
from the \thmref{dixon} it follows that $g$ is a isomorphism, so proving the claim.

It is a natural question if the morphism $p\colon R_n\to M_n$ (see (\ref{defw})) is injective, posed  by
Ellinsgrud and Stromme in \cite{ES} and settled again in \cite{PT}.
This question is equivalent to ask how many complete  $(n+1)$-gons
are inscribed in the general Darboux curve of degree $n$.

Now observe that
the secant order $d_{n+1}(X)$ (see def. \ref{secord}) can be computed as degree of the induced  morphism
on quotients
$$R_n=\sigma^{n+1}(X)//SL(V)\rig{}\sigma_{n+1}(X)//SL(V)$$

We have the factorization

$$\begin{array} {ccc}
R_n\\
\dow{\pi} &\sear{p}\\
\sigma_{n+1}(X)//SL(V)&\rig{h}&M\\
\end{array}$$

and it follows that

$$\deg p = d_{n+1}(X) \cdot \deg h$$

$\deg h$ measures how many different theta divisors defined by
inscribed $(n+1)$-gons (like in \propref{nn+1} lie over the generic Darboux curve.

The degree of $h$ is not known, although it is expected that it is one.

The secant order $d_{n+1}(X)$ measures how many $(n+1)$-gons give linearly equivalent $\theta$ divisors
like in \propref{nn+1}.

The following reformulation of the \thmref{maincc} in this setting looks interesting.

\begin{thm0}\label{reform} Let $n\ge 5$. $d_{n+1}(X)=1$  unless
$X= \P2\times\P{{n-1}}$
embedded with the linear system $\O (1,2)$ is $(n+1)$-weakly defective.
\end{thm0}
\proof  By \thmref{maincc} and \thmref{expn+1}.
\qedd

We do not know if $X$ in the above theorem is $(n+1)$-weakly defective,
the \thmref{expn+1} says only that $X$ is not $(n+1)$-defective.

M. Toma \cite{Tom} has shown a uniqueness result for curves such that their inscribed $(n+1)$-gon is also circumscribed to a 
smooth conic (and $n\ge 5$), they are called Poncelet curves.

\section{The moduli space of symplectic vector bundles on the plane}

Let $M(r,n)$ be the moduli space of stable bundles of rank $r\ge 2$ on $\P 2$
with $(c_1,c_2)=(0,n)$.
It is known that
$M(r,n)$ is not empty if and only if $ r\le n$.
$M(r,n)$ is a smooth irreducible variety of dimension $2rn-r^2+1$ (\cite{Hu}).

Moreover for the general $E\in M(2,4)$  we have $H^0(E(1)=2$,
  while for  $n\ge 5$ the Brill-Noether locus
 $H_n=\{E\in M(2,n) | h^0(E(1))\ge 1\}\subseteq M(2,n)$ is a irreducible
 subvariety of dimension $3n+2$, which is proper if $n\ge 6$ (\cite{Bar}). The bundles $E\in H_n$ are called Hulsbergen bundles. 
 
 In this and in next section we will reprove the main results of \cite{Bar} and \cite{Hu}
in the case of symplectic bundles, by connecting them with the higher secant varieties.

\begin{defn0} A vector bundle is called symplectic if there is a isomorphism
$$\alpha\colon E\to E^{\vee}$$
such that $\alpha^t=-\alpha$. 
\end{defn0}
In particular it follows that if $E$ is symplectic then $c_1(E)=0$, $r$ is even and
$\wedge^2E$ contains $\O$ as a direct summand.
If moreover $E$ is stable we have that $h^0(\wedge^2 E)=\C$, so that the isomorphism
$\alpha$ is unique up to scalar multiplication.

Let now $E$ be a stable vector bundle of rank $r$ on $\P 2$ with $c_1(E)=0$, $c_2(E)=n$.
A simple computation shows that $\chi(E(-1))=\chi(E(-2))=-n$ is independent by $r$.
 
 \begin{thm0}
 \label{barthmonad}
 Let $E$ be a symplectic bundle of rank $r$ on $\P 2=\P{{}}(U)$ with $c_2(E)=n$
 such that $H^0(E)=0$.
 Denote $V=H^1(\P 2,E(-1))$ which is a vector space of dimension $n$.
 Then $E$ is the cohomology bundle of the following Barth monad
$$I\otimes \O\rig{g}V^{\vee}\otimes\Omega^1(2)\rig{f}V\otimes \O (1) $$
that is $E=\textrm{Ker\ } f/\textrm{ Im\ } g$
where 
 $f\in U\otimes S^2V$ is the natural (symmetric) multiplication map
and $I=H^1(E(-3))=H^1(E)^{\vee}$ has dimension $r-n$.

Conversely the cohomology bundle $E(f)$ of such a monad 
where $f\in U\otimes S^2V$ is
 a symplectic bundle of rank $r$ with $c_2(E)=n$
 such that $H^0(E)=0$.
\end{thm0}

\proof
The Beilinson table for $E(-1)$ is

$$\begin{array}{lll}0&0&0\\
I&V^{\vee}&V\\
0&0&0\\
\end{array}$$
hence by twisting by $\O(1)$ we get from \thmref{beilinson} the monad in the statement.
Note that $f\in Hom(V^{\vee}\otimes\Omega^1(2), V\otimes \O (1))=U\otimes V\otimes V$.

Since Serre duality is induced by cup product which is skew commutative in odd dimension,
it is well known that $f=f^t$ (for details see \cite{Bar} Prop. 1),
hence $f\in U\otimes S^2V$. The converse is trivial.\qedd

\begin{prop0}\label{slv}
Two simple bundles $E(f)$, $E(f')$ as
in \thmref{barthmonad} are isomorphic if and only if
$f$, $f'$ are $SL(V)$-equivalent.
\end{prop0}
\proof It is easy to check that a morphism between bundles lifts to a morphism between the corresponding Barth monads.
The details are left to the reader. \qedd
From the Barth monad of \thmref{barthmonad} we get the cohomology map
$$H^0(f)\colon V^{\vee}\otimes H^0(\Omega^1(2))\to V\otimes H^0( \O (1))$$

It is important to remark that $H^0(f)$ 
is identified with $S_f\colon U\otimes V^{\vee}\to U^{\vee}\otimes V$ in formula (\ref{sfi}) of section $4$.
Note that $\textrm{rank\ }H^0(f)=2n+r$, indeed the kernel of $H^0(f)$ contains $I$ which has dimension $r-n$,
and it cannot be bigger,
otherwise $H^0(E)\neq 0$ contradicting the stability of $E$.

We denote $M_{sp}(r,c_2)\subseteq M(r,c_2)$ the moduli space of stable symplectic vector bundles
(note that $r$ is even).
We recall that the adjoint representation for the symplectic group $Sp({\bf C}^r)$
is isomorphic to the symmetric power $S^2{\bf C}^r$.

Since any $E\in M_{sp}(r,c_2)$ is simple, we get
$h^0(S^2E)=0$, moreover $h^2(S^2E)=h^0(S^2E(-3))=0$.

We have
$ c_2(S^2E)=n(r+2) $
and  we get by Hirzebruch-Riemann-Roch theorem

$$h^1(S^2E)= -\chi (S^2E)=(r+2)n-{{r+1}\choose 2}$$

It follows that $M_{sp}(r,n)$ (when nonempty) is smooth of dimension
$(r+2)n-{{r+1}\choose 2}$.

Following \cite{Hu}, denote $M_{sp}^0(r,n)=\{E\in M_{sp}(r,n) | E_{|l}=\O^r \textrm{for some line\ }l \}$.
We remark that, by semicontinuity, if $E_{|l}$ is trivial on a line $l$,
then it is trivial on the general line $l$.

\begin{prop0}
\label{hirschowitz}

$$\overline {M_{sp}^0(r,n)}= M_{sp}(r,n)$$
\end{prop0}

\proof (Hirschowitz) Let $E\in M_{sp}(2h,n)$. For a line $l$ we have a small deformation
of $E_{|l}$ such that $E_{|l}$ is trivial. We have the exact sequence
$$0\rig{}S^2E(-1)\rig{}S^2E\rig{}S^2E_{|l}\rig{}0$$
which yields
$$H^1(S^2E)\rig{}H^1(S^2E_{|l})\rig{}0$$

Hence the deformation on $l$ lifts to a deformation on ${\bf P}^2$.\qedd

Let $E=E(f)$ like in the statement of \thmref{barthmonad}.
The exact sequence

$$0\rig{}E(-2)\rig{}E(-1)\rig{}E(-1)_{|l}\rig{}0$$
yields

$$0\rig{}H^0(E(-1)_{|l})\rig{}H^1(E(-2))\rig{}H^1(E(-1))$$

hence $E_{|l}$ is trivial if and only if 
the discriminant of the morphism $H^1(E(-2))\rig{}H^1(E(-1))$
is nonzero, that is if and only if $\Delta(f)$ evaluated at $l$
is nonzero.

Hence $E(f)$ is trivial on the general line $l$ if and only if 
$\Delta(f)$ is not identically zero, that is if and only if $f$ corresponds
to a semistable point in $\P{{}}(U\otimes S^2V)$ for the $SL(U)\times SL(V)$-action.

\begin{defn0}
$$K_{r,n}=\{f\in P{{}}(U\otimes S^2V) | \textrm{rank\ }H^0(f)=2n+r\}$$
\end{defn0}
The general element in $\sigma_{n+(r/2)}(X)$
belongs to $K_r$ (so it is nonempty).

 In order to prove the irreducibility of $M_{sp}(r,n)$ we need the following auxiliary result
 from \cite{BPV} cor. 3.6 
(see also \cite{Bas} cor. 2.6 for a elementary proof in the spirit of \cite{Hu}).
 
\begin{prop0}[Brennan-Pinto-Vasconcelos] 
\label{bpv}
Let $V$ be a vector space of dimension $n$
and let $r$ be even.
The subvariety $$J_{r,n}:=\{(P,Q)\in S^2V\times S^2V | \textrm{rank\ }[P,Q]=r\}$$
is irreducible of codimension
${n-r}\choose{2}$, and moreover its reduced equations are the pfaffians of order $r+2$
of $[P,Q]$.
\end{prop0}

\begin{thm0}
The moduli space $M_{sp}(r,n)$ is  irreducible of dimension $(r+2)n-{{r+1}\choose 2}$.
\end{thm0}

\proof We have the open subvariety $\tilde K_{r,n}\subseteq K_{r,n}$ defined by $f\in K_{r,n}$
such that the corresponding morphism
$$V^{\vee}\otimes\Omega^1(2)\rig{f}V\otimes \O (1) $$
is surjective. Any such $f$ defines $E(f)$ as cohomology bundle of the corresponding Barth monad.
It is easy to see that there is a universal bundle ${\cal E}$ over $\P{{}}(U)\times \tilde K_{r,n}$
such that ${\cal E}_{{\bf P}^2\times\{ f\} } = E(f)$, by constructing a universal monad like in \cite{Hu}
prop. 1.6.1. Since stability is a open property, we have
a open subvariety $\tilde K_{r,n}^s\subseteq \tilde K_{r,n}$ consisting of $f$ such that $E(f)$ is stable.
By the universal property of moduli space we have a surjective morphism
$\tilde K_{r,n}^s\rig{\pi} M_{sp}(r,n)$.
By \propref{hirschowitz} it is enough to prove that
$M_{sp}^0(r,n)$ it is irreducible, hence it is enough to prove that 
$\pi^{-1}(M_{sp}^0(r,n))=\tilde K_{r,n}^s\setminus Z(\Delta)$ is irreducible.
We will prove that $\tilde K_{r,n}\setminus Z(\Delta)$ is irreducible.
For any $x\in\P{{}}(U)$ define $$\tilde K_{r,n,x}=\{ f\in K_{r,n} | \Delta (f)(x)\neq 0\}$$

These are open subsets in $\tilde K_{r,n}$ such that, for any $x, y\in\P{{}}(U)$,
$\tilde K_{r,n,x}\cap\tilde K_{r,n,y}$ is a non empty open subsvariety of $\tilde K_{r,n}$
and moreover
$$\bigcup_{x\in\P{{}}(U)}\tilde K_{r,n,x}=\tilde K_{r,n}\setminus Z(\Delta)$$

Finally we prove (thanks to the $SL(U)$-action) that for $z=(0,1,0)$ $\tilde K_{r,n,z}$ is irreducible.
Let $C\tilde K_{r,n,z}$ be the affine cone over $\tilde K_{r,n,z}$.
In the matrix representation of \lemref{eqz} this means that the slice called $Q$ of $S_f$ is invertible,
hence we have a fibration $C\tilde K_{r,n,z}\to S(n)$ sending $f$ to $Q$,
where $S(n)$ is the variety of symmetric invertible matrices of order $n$.
This fibration is $GL(V)$-equivariant and all its fibers are isomorphic to $J_{r,n}$ defined
in \propref{bpv}. Hence the result follows by \propref{bpv}.
Note that 
$$\dim S(n)+\dim J_{r,n}-\dim GL(V)=3{{n+1}\choose 2}-{{n-r}\choose{2}}-n^2= (r+2)n-{{r+1}\choose 2}$$
\qedd

{\bf Question} Are the moduli spaces of orthogonal bundles on $\P 2$ irreducible ?

\section{The Barth map and Brill-Noether loci}

We keep all the notations from the previous section.

Let remark the following result of LePotier

\begin{prop0}[LePotier]Assume $f\in\P{{}}(U\otimes S^2V)$ is a semistable point (see \propref{wall})
for the $SL(V)$-action.
If $E(f)$ is a stable bundle, then $f$ is a stable point.
\end{prop0}

\proof The proof is essentially the same as in \cite{P2} Lemma 2, 
and appears in different forms also in \cite{Bar} and \cite{Hu}, so we omit it. \qedd

Now we consider the GIT quotient $\tilde K_{r,n}^s//SL(V)$. Since all points
not lying in $Z(\Delta)$ are $SL(V)$-semistable, by \propref{slv} and the above construction
we have a surjective morphism $\tilde K_{r,n}^s//SL(V)\rig{\tau} M_{sp}^0(r,n)$
which is birational. Consider the closure

$$\overline{K_{r,n}}=\{f\in P{{}}(U\otimes S^2V) | \textrm{rank\ }H^0(f)\le 2n+r\}=\bigcup_{s\le r}K_{s,n}$$

We get the projective scheme
$$M_{sp}^{mon}(r,n):=\overline{K_{r,n}}//SL(V)$$
where the suffix $mon$ reads for monads, which has the following properties

(i) it is birational to the Maruyama moduli space $\overline{M_{sp}(r,n)}$

(ii) the morphism $\Delta$ factors through
$$\begin{array} {ccc}
\overline{K_{r,n}}\setminus Z(\Delta)\\
\dow{\pi} &\sear{\Delta}\\
M_{sp}^{mon}(r,n)&\rig{b_{r,n}}&\P{{}}(S^nU)\\
\end{array}$$ 

(iii) The following diagram commutes
$$\begin{array} {ccc}
\tilde K_{r,n}^s//SL(V)&\rig{\tau} &M_{sp}^0(r,n)\\
\dow{i} &&\dow{J}\\
M_{sp}^{mon}(r,n)&\rig{b_{r,n}}&\P{{}}(S^nU)\\
\end{array}$$ 
where $i$ is an open embedding and $J(E)$ is the degree $n$ curve of jumping lines of $E$ supported by
$$\{ l\in P(U^{\vee})| E_{|l}\neq\O^r\}$$

\begin{defn0}
The morphism $b_{r,n}\colon M_{sp}^{mon}(r,n)\rig{}\P{{}}(S^nU)$
which factors in (ii) above is called the Barth map.
\end{defn0}

Its degree and the degree of its image are equal to the ones for Barth maps
as defined in \cite{PT} in the case $r=2$. Note that from property (ii) above
the Barth map can be computed by the symmetric determinantal morphism $\Delta$.
See \cite{P2}
for the connection with Donaldson invariants.

\begin{prop0}
The Barth map is generically finite, the image of the Barth map is irreducible
and has codimension
$1+\frac{(n+1-r)(n-2-r)}{2}={{n-r}\choose{2}}$ in the projective space of plane curves of degree $n$.
\end{prop0}

\proof By \thmref{dixon} and by (ii) above. \qedd

\begin{coro0}
The Barth maps $b_{n,n}$ ($n$ even) and $b_{n-1,n}$ ($n$ odd) are dominant.
\end{coro0}

The case $r=n-2$, $n$ even looks particularly interesting because
the image of $b_{n-2,n}$ is a hypersurface in $M_n$.
For $n=4$ this hypersurface is the L\"uroth hypersurface of section 6.
By \thmref{simstrassen} we get that , with $n$ even and $k=\frac{3n}{2}-1$
\begin{equation}
\label{n-2}
\overline{K_{n-2,n}}=\sigma_{k}(X)
\end{equation}

We remark that $\sigma_{n+1}(X) = \overline{K_{2,n}}$ for $n=4, 5$,
while for $n\ge 6$ we have $\dim \overline{K_{2,n}}=\dim \sigma_{n+1}(X) +(n-5)$.
Note also that note that  
$\sigma_{n-1}(X)\subseteq Z(\Delta)$ but this inclusion is strict for $n\ge 3$.

Consider the diagram
$$\begin{array} {ccc}
M_{sp}^{mon}(r,n)\\
\dow{i} &\sear{b_{r,n}}\\
\P{{}}(U\otimes S^2V)//SL(V)&\rig{g}&\P{{}}(H^0(\P2,\O(n)))\\
\end{array}$$

By \propref{wall} and by \thmref{dixon} it follows that 
$g$ is generically finite, in particular the fiber over a smooth $C\in M_n^{sm}$
is given by the theta characteristic $\{L\in Pic(C)| L^2=K_C, h^0(L)=0\}$ .

\begin{prop0}[Beauville, Catanese]\label{exactly}
Let be given a general smooth plane curve $C$ 
and genus $g={{n-1}\choose {2}}$. Then the set
$\{L\in Pic(C)| L^2=K_C, h^0(L)=0\}$ has cardinality
$$\left\{\begin{array}{cl}2^{n-1}(2^n+1)&\textrm{if $n$ is even or if $n\equiv 3,5$\ mod\ $8$ }\\
2^{n-1}(2^n+1)-1&\textrm{if $n\equiv 1,7$\ mod\ $8$}\\ 
\end{array}\right.$$
\end{prop0}
\proof Let $n$ be even. The moduli space ${\cal T}_n$ of pairs $(C,\theta)$ where
$C$ is a smooth plane curve of degree $n$ and $\theta$ is a theta characteristic
has exactly two irreducible components ${\cal T}_n^0$ and ${\cal T}_n^1$,
 see \cite{Bea2} prop. 3, corresponding to the parity of $\theta$, which are both
\`etale covering of the space $U_n\subset \P{{}}(H^0(\P 2,\O (n)))$ of smooth plane curves.

By the \thmref{dixon}, the generic plane curve $C$ has a $\theta$ such that $h^0(\theta)=0$.
It follows that the subvariety in  ${\cal T}_n^0$  of pairs  $(C,\theta)$
such that $h^0(\theta)\ge 2$ has codimension at least one in ${\cal T}_n^0$,
and the same is true for its projection on $U_n$. The number of sheets of
${\cal T}_n^0$ has been computed in \cite{At}.
For $n$ odd ${\cal T}_n$ has a third irreducible component corresponding to
$\O(\frac{n-3}{2})$ which has always $h^0\ge 2$ and it is even if 
$n\equiv 1,7$\ mod\ $8$. \qedd

It is well known that, in the moduli space $M_g$ of curves of genus $g$, 
the locus of curves which have a even theta characteristic such that
$h^0(\theta)\ge 2$ has pure codimension one, see \cite{Teix} theor. 2.16.
This is called the {\it theta locus}.

For $g=3$ this divisor coincides with the hyperelliptic locus and so it does not meet
the locus of plane curves of degree $4$.

For plane curves of even degree $\ge 6$ the situation changes, thanks to the following
interesting examples.

\begin{prop0}[Pirola]\label{pirola}
Let $h\ge 3$. Consider three general curves of degree respectively $h$, $2h-3$, $3$ with equation
respectively
$C_h$, $C_{2h-3}$, $C_3$, $h\ge 3$.
Then the curve $K$ with equation
$C_h^2-C_{2h-3}C_3$ 

(i) is smooth  of degree $2h$

(ii) it has a theta characteristic $\theta$ such that $h^0(\theta)=1+{{h-1}\choose{2}}$,
in particular $\theta$ is even for $h=0,3\ \textrm{mod\ }4$
and $\theta$ is odd for $h=1,2\ \textrm{mod\ }4$.
\end{prop0}
\proof 
Take $C_h=x_0^h$ and factor
$x_1^{2h}+x_2^{2h}=F_{2h-3}(x_1,x_2)F_3(x_1,x_2)$.

Posing  $C_{2h-3}=F_{2h-3}$, $C_3=F_3$, we get that $K$ is the Fermat curve,
hence it remains smooth by deforming the three equations, so proving (i).

In order to prove (ii) we denote by $\theta$ the divisor on $K$ given by $\{C_h=0\}\cap \{C_{2h-3}=0\}$ and 
we observe that $K\cap\{C_{2h-3}=0\}$
is contained in $K\cap \{C_h=0\}$. Hence the curve with equation $C_{2h-3}$ cuts $K$ in the divisor $2\theta$,
which means that $\theta$ is a theta characteristic. By Riemann Roch
$h^0(\theta)=h^0((2h-3)H-\theta)$
and from the sequence
$$0\rig{}\O_{{\bf P}^2}(-3)\rig{}{\cal I}_{\theta,{\bf P}^2 } (2h-3)\rig{}
\O_K(  (2h-3)H-\theta ) \rig{}0$$

we compute $h^0(K,(2h-3)H-\theta)=h^0({\cal I}_{\theta,{\bf P}^2 } (2h-3))$
which in turn can be computed by the Koszul complex
                    
$$ 0\rig{}\O_{{\bf P}^2}(-h) \rig{}\O_{{\bf P}^2}(h-3)\oplus\O_{{\bf P}^2}\rig{}{\cal I}_{\theta,{\bf P}^2 } (2h-3)\rig{}0$$

and it is $1+{{h-1}\choose{2}}$ as we wanted. \qedd

These examples by Pirola show that the theta locus actually meets
 the variety of smooth plane curves of degree $n$ for any $n\ge 5$.
 
 Also the even theta locus meets the variety of smooth plane curves of even degree $2h$ ($h\ge 3$),
with the possible exceptions $h=1,2\ \textrm{mod\ }4$.

\begin{prop0}\label{thetaquartic}
Every smooth plane quartic has exactly $36$ theta characteristic
such that $h^0(\theta)=0$ and $28$ theta characteristic
such that $h^0(\theta)=1$
\end{prop0}

\proof This is well known and follows from \cite{At}
and the remark that, by Clifford's theorem, on a smooth plane quartic every theta-characteristic satisfies
$h^0(\theta)\le 1$.\qedd

 The degree of $b_{2,n}$
is one by  \cite{PT}.
The degree of the image of $b_{2,n}$ are known, see e.g. \cite{EG}.

The {\it proper transform} $\Delta^*(W)$ of a subvariety $W\subseteq M_n$
through the morphism
$$\Delta\colon\P{{}}(U\otimes S^2V)\setminus Z(\Delta)\rig{} M_n$$
is defined as the Zariski closure
of $\Delta^{-1}(W)\setminus Z(\Delta)$.
The degree of the proper transform is difficult to compute,
with the exception of hypersurfaces. Assume $W$ is a hypersurface.
Since $Z(\Delta)$ has bigger codimension, we get that $\Delta^*(W)$
is a hypersurface too and
$$\deg \Delta^*(W)= n \deg W$$

\begin{thm0}\label{ddelta} 
Let $E$ be a symplectic vector bundle on ${\bf P}^2$
and let $r=n-2$, hence the image of the Barth map is a hypersurface.
Let $g={{n-1}\choose 2}$.
Two mutually exclusive cases are possible.

(i) The image of the Barth map $b_{n,n-2}$ is contained in the theta locus

(ii) $\textrm{(degree of\ } b_{n,n-2})\cdot\textrm{(degree of image of\ } b_{n,n-2}) =3\cdot 2^{g-2}(2^g+1)$
\end{thm0}

\proof Assume that (i) does not hold. Then the generic curve in the image of the Barth map
has $2^{g-1}(2^g+1)$ theta-characteristic $\theta$ such that $h^0(\theta)=0$.

Denote $a=\textrm{(degree of\ } b_{n,n-2})$, $b=\textrm{(degree of image of\ } b_{n,n-2})$.
The proper transform of the image of Barth map is a hypersurface in $\P{{}}(U\otimes S^2V)$
of degree $nb$ which contains $\sigma_{\frac{3n}{2}}(X)$ which has degree $\frac{3n}{2}$
by the remark at the end of section $4$.

Consider the intersection with the proper transform of a general line.
It is given by $b$ fibers (over smooth curves), each of them is the union of $a$ $SL(V)$-orbits,
which have the same degree.
Hence we have the equation
$$\frac{nb}{3n/2}=\frac{2^{g-1}(2^g+1)}{a}$$
which gives the thesis.
\qedd

We do not know which of the two possibilities hold, except for $n=4$ where we have the following
corollary, that was proved first in \cite{PT} for any $c_2$. Our approach is different.

\begin{coro0}
The Barth map $b_{2,4}$ is generically injective.
\end{coro0}
\proof We know that by \propref{thetaquartic} the case (i) of \thmref{ddelta} cannot occur.
By \cite{Mo} we know that the degree of the L\"uroth hypersurface, which is the image of $b_{2,4}$, is $54$.
The result follows. \qedd

\begin{rem0} Although we do not know the explicit expression of the L\"uroth invariant $L$ of degree $54$, we can say
that its pullback $\Delta^*L$ has degree $216$ and it contains $Pf(S_{f})$ of degree $6$
as irreducible factor.
\end{rem0}
\begin{thm0}\label{jumpn-2} Let $E\in M_{sp}(n-2,n)$ be general and let $C=J(E)$ be its curve of jumping lines,
so that $C$ lies in the image of $b_{n-2,n}$. Then there are  $r_1,\ldots ,r_k$ lines where $k=\frac{3n}{2}-1$
and  linear forms $h_i$ such that the equation of $C$ can be written as
$$\Delta(\sum_{i=1}^kr_ih_i^2)=0$$
Moreover the varieties of  lines $r_1,\ldots ,r_t$ which describe $C$ in the above equation has dimension $\frac{n}{2}-1$.
\end{thm0}

\proof By (\ref{n-2}) we find the $k$ lines and the $k$ linear forms.
By \thmref{simstrassen} and \propref{infty} we compute
$$(\frac{3n}{2}-1)(n+2)-1-\left[3{{n+1}\choose{2}}-2\right]= \frac{n}{2}-1 $$  \qedd

 When $n=4$ the \thmref{jumpn-2} reduces to the result of \cite{Bar}
that the jumping lines of a general $E\in M_{sp}(2,4)$ are L\"uroth quartics.

The linear forms $h_i$ define a $k\times n$ matrix $H$.
Let $I$ be any subset of $n$ rows and let $h_I$ be the corresponding minor of $H$.

Then the equation of  $C$ can be written as
$$\sum_{I}h_I^2\prod_{j\in I}r_j=0$$

{\bf Question} What is the geometric interpretation of the curves of degree $n$ lying in the image of the Barth map
$b_{n-2,n}$ ? For $n=4$ they are the L\"uroth quartics. 
How is $\theta$ related to the data of the equation ?
For $r=3$ they are sextics with a determinantal representation
arising from $9$ lines. 

\medskip

The last case we are interested regards the bundles
defined from a general $f\in\sigma_{n+(r/2)}(X)$.

Let $E=E(f)\in M_{sp}(r,n)$ on $\P{{}}(U)$. We remark that from the Barth monad we have

$$h^0(E(1))=\dim\ker \left[V^{\vee}\otimes H^0(\Omega^1(3))\rig{H^0(f(1))}V\otimes H^0(\O (2))\right]-3(n-r)$$

It is well known that $\Omega^1(3)$ is the tangent bundle, then
$H^0(\Omega^1(3))=\textrm{ad\ } U^{\vee}$, moreover $H^0(\O (2))=S^2 U^{\vee}$.

\begin{lemma0}
\label{4term}
On $\P{{}}(U)$ we have the exact $SL(U)$-homogeneous sequence
$$0\rig{}\O(-2)\rig{}S^2U\rig{q}\textrm{ad\ }U(1)\rig{}U(2)\rig{}0$$
\end{lemma0}
\proof It
is a particular case of the Four Term Lemma (Lemma 26 in \cite{OR}),
it can also be found by an explicit computation.\qedd

\begin{thm0}\label{hulsbergen}
Let $f\in\sigma_{n+(r/2)}(X)$ general and $E=E(f)\in M_{sp}(r,n)$ on $\P{{}}(U)$.
Then $h^0(E(1))\ge r/2$.
\end{thm0}

\proof Consider first $f=u\otimes v^2\in X$. The $n\times n$ matrix  corresponding
to the map $V^{\vee}\otimes \Omega^1(3)\rig{}V\otimes \O (2)$ has only one non zero coefficient,
which is $u$, say at the entry $(1,1)$. At level of $H^0$, the contraction by $u$ corresponds
at the evaluation of $q^t$ of \lemref{4term} at $u$.
From \lemref{4term} it follows that $\textrm{codim\ }\ker H^0(f(1))=5$.

If $f=\sum_{i=1}^{n+(r/2)}f_i$ with $f_i\in X$ we get that

$$ \cap_i \ker H^0(f_i(1))\subseteq \ker H^0(f(1))$$

hence $$\textrm{codim\ }\ker H^0(f(1))\le \sum_i \textrm{codim\ }\ker H^0(f_i(1))=5\left[n+(r/2)\right]$$

It follows that $$h^0(E(1))\ge 8n-5\left[n+(r/2)\right]-3(n-r)=r/2$$
as we wanted. \qedd

\begin{rem0} The \thmref{hulsbergen}  is meaningful only when $n>\frac{5r}{2}$,
otherwise all $E\in M_{sp}(r,n)$ satisfy the inequality $h^0(E(1))\ge r/2$,
because $\chi(E(1))=3r-n$.
\end{rem0}
\begin{rem0} When $r=2$ the bundles constructed in the \thmref{hulsbergen}
are exactly the Hulsbergen bundles in \cite{Bar}. Their curve of jumping lines is a Darboux curve
by \propref{darblocus}. When $n=4$ the intersection computed in the above proof is not transversal
(the reader can recognize here the flavour of the Terracini lemma), indeed
in such a case the general $E$ satisfies $h^0(E(1))=2$.
It is easy to show that their section vanishes exactly on the vertices
of the inscribed $(n+1)$-gon.
\end{rem0}
Let $M_{sp}(r,n)^k:=\{ E\in M_{sp}(r,n)|h^0(E(1)\ge k\}$ be the Brill-Noether locus and consider $E\in M_{sp}(r,n)^k$.
The Brill-Noether theory says that the tangent space to $M_{sp}(r,n)^k$ at $E$ is the kernel
of the natural morphism
$$H^1(S^2E)\rig{}H^0(E(1))^{\vee}\otimes H^1(E(1))$$

If $h^0(E(1))=r/2$, such tangent space has dimension bigger or equal than
$$ (r+2)n-{{r+1}\choose 2}-(r/2)\left[n-(5r/2)\right]=n\left[2+(r/2)\right]+\frac{3r^2-2r}{4}$$
which is equal to $$\dim\Delta(\sigma_{n+(r/2)}(X))=n\left[2+(r/2)\right]+r$$
only when $r=2$, otherwise it is bigger. This means that the general bundle in 
$M_{sp}(r,n)^{(r/2)}$ comes from $\sigma_{n+(r/2)}(X)$ only if $r=2$,
which is the case of Hulsbergen bundles.
  
\bigskip
{\sf

Giorgio Ottaviani\\
Dipartimento di Matematica U. Dini, Universit\`a di Firenze\\
viale Morgagni 67/A,  50134 Firenze, Italy\\
ottavian@math.unifi.it
}


\bigskip
{\small


\begin{thebibliography}{Dilloo Dilloo 33}

\bibitem[AOP]{AOP} H. Abo, G. Ottaviani, C. Peterson,  {\em Induction for secant varieties of Segre varieties},  math.AG/0607191

\bibitem[AO]{AO}  V. Ancona, G. Ottaviani, Some applications of Beilinson theorem to projective
spaces and quadrics , Forum Math. 3 (1991), 157-176

\bibitem[At]{At} M. Atiyah, Riemann surfaces and spin structures. 
Ann. Sci. École Norm. Sup. (4) 4 (1971) 47-62

\bibitem[Bar]{Bar} W. Barth, Moduli of vector bundles on the projective plane , Invent. Math. 42 (1977),
63-92

\bibitem[Bas]{Bas} R. Basili, On the irreducibility of varieties of commuting matrices, 
J. Pure Appl. Algebra, 149 (2000) 107-120

\bibitem[Bea1]{Bea1} A. Beauville, Determinantal hypersurfaces, 
Michigan Math. J. 48 (2000), 39-64 

\bibitem[Bea2]{Bea2} A. Beauville, Le groupe de monodromie des familles universelles d'hypersurfaces at 
d'intersections compl\`etes, Complex analysis and Algebraic Geometry, LNM 1194, 195-207, Springer 1986

\bibitem[B]{B} A. Beilinson, Coherent sheaves on $\P n$ and problems of linear algebra,
 Funkt. Analiz Prilozhenia, 12 n.3, 68-69 (1978)

\bibitem[BPV]{BPV} J.P. Brennan, M.V. Pinto, W.V. Vasconcelos,
The Jacobian module of a Lie algebra, 
Trans. Amer. Math. Soc. 321 (1990), no. 1, 183-196. 

\bibitem[CarCa]{CarCa} E. Carlini, M.V. Catalisano, Existence results for rational normal curves , math.AG/0603137


\bibitem[CaCh]{CaCh} E. Carlini, J.V.Chipalkatti,
 On Waring's problem for several algebraic forms, Comment. Math. Helv. 78 (2003), no. 3, 494-517.

\bibitem[CGG1]{CGG1} M.V. Catalisano, A.V. Geramita, A. Gimigliano,  {\em Ranks of
tensors, secant varieties of Segre varieties and fat points}, Linear Algebra
Appl. 355 (2002), 263-285.  Erratum: Linear Algebra Appl. 367 (2003), 347-348.

\bibitem[CGG2]{CGG2} M.V. Catalisano, A.V. Geramita, A. Gimigliano,  {\em Secant
varieties of Segre-Veronese varieties}, Projective varieties with unexpected properties, 81-107,
 W. de Gruyter , Berlin, 2005.

\bibitem[CC1]{CC1} L. Chiantini, C. Ciliberto, {\em Weakly defective varieties}, Trans.
AMS 354 (2002), 151-178.

\bibitem[CC2]{CC2} L. Chiantini, C. Ciliberto,  {\em On the concept of $k$-secant order of a variety},
 J. London Math. Soc. (2) 73 (2006), no. 2, 436-454.

\bibitem[Ci]{Ci} C. Ciliberto, {\em Geometric aspects of polynomial interpolation in
more variables and of Waring's problem},
ECM Barcelona 2000, vol. I, Progr. Math. 201 (2001). 289-316

\bibitem[D]{D} I. Dolgachev, Topics in Classical Algebraic Geometry, lecture notes

\bibitem[DK]{DK} I. Dolgachev, V. Kanev, Polar covariants of plane cubics and quartics. 
Adv. Math. 98 (1993), no. 2, 216-301

\bibitem[ES]{ES} G. Ellingsrud, A. Stromme, Bott's formula and enumerative geometry,
J. Amer. Math. Soc. 9 (1996), no. 1, 175--193. 

\bibitem[EG]{EG} G. Ellingsrud, L.  G\"ottsche, 
Variation of moduli spaces and Donaldson invariants under change of polarization,
J. Reine Angew. Math. 467 (1995), 1-49

\bibitem[FGA]{FGA} A. Grothendieck, Fondements de la G\'eometrie Alg\'ebrique, S\'eminaire Bourbaki 1957-62,
Secr\'etariat Math., Paris (1962)

\bibitem[Fr]{Fr} W. Frahm, Bemerkung \"uber das Fl\"achennetz zweiter Ordnung, Math. Ann. 7, 635-638 (1874)

\bibitem [Hu]{Hu} K. Hulek , On the classification of stable rank-$r$ vector bundles over the projective plane.
 Vector bundles and differential equations (Proc. Conf., Nice, 1979), pp. 113-144, Progr. Math., 7, Birkh\"auser,
 Boston, 1980

\bibitem[LM]{LM}  J.M. Landsberg, L. Manivel, Generalizations of Strassen's equations for secant
 varieties of Segre varieties, AG/0601097

\bibitem[LW]{LW} J.M. Landsberg, J. Weyman, On the ideals and singularities of secant varieties of Segre varieties
, math.AG/0601452

\bibitem[Mo]{Mo} F. Morley, On the L\"uroth quartic curve, Amer. J. Math. 36 (1918), 279-2820 

\bibitem[OR]{OR} G. Ottaviani, E. Rubei, Resolutions of homogeneous bundles on ${\bf P}^2$, math.AG/0401405,
Annales de l'Institut Fourier 55 (2005),973-1015 

\bibitem[OSS]{OSS} C. Okonek, M. Schneider, H. Spindler, Vector bundles on complex projective spaces,
Birkhauser, Boston 1980


\bibitem[P1]{P1} J. Le Potier, Sur l'espace de modules de fibr\'es de Yang et Mills, , S\'em ENS 1979-82,
Progress in Math. 37, Birkh\"auser, Boston 1983

\bibitem[P2]{P2} J. Le Potier, Fibr\'e d\'eterminant et courbes de saut sur les surfaces alg\'ebriques,
Complex projective geometry (Trieste, 1989/Bergen, 1989), 213-240, 
London Math. Soc. Lecture Note Ser., 179, Cambridge Univ. Press, Cambridge, 1992 

\bibitem[PT]{PT} J. Le Potier, A. Tikhomirov, Sur le morphisme de Barth, 
Ann. Sci. \'Ecole Norm. Sup. (4) 34 (2001), no. 4, 573-629 

\bibitem[Ru]{Ru} F. Russo, {\em Tangents and Secants to Algebraic Varieties},
Publicacoes Matematicas do IMPA. $24$ Colloquio Brasileiro de Matematica. (IMPA),
Rio de Janeiro, 2003.

\bibitem[S]{S} V. Strassen, Rank and optimal computation of generic tensors, Linear Algebra Appl. 52 (1983) 645-685. 

\bibitem[Teix]{Teix} M. Teixidor i Bigas, Half-canonical series on algebraic curves,
Trans. Amer. Math. Soc. 302 (1987), no. 1, 99-115

\bibitem[T]{T} E. Toeplitz , \"Uber ein Fl\"achennetz zweiter Ordnung, Math. Ann. 11 (1877),  434-463

\bibitem[Tom]{Tom} M. Toma, Birational models for varieties of Poncelet curves, 
Manuscripta Math. 90 (1996), no. 1, 105-119 

\bibitem[W]{W} C.T.C. Wall, Nets of quadrics and theta-characteristics of singular curves,
Philos. Trans. Roy. Soc. London Ser. A 289 (1978), no. 1357, 229-269

\bibitem[Z]{Z} F.L. Zak, {\em Tangents and Secants of Algebraic Varieties},
Translations of Mathematical Monographs, 127.
AMS, Providence, RI, 1993.



\end{thebibliography}}
\end{document}